\documentclass{article}
\usepackage{hyperref}
\usepackage{amsthm}
\usepackage{amsmath,amssymb}
\newtheorem{theorem}{Theorem}[section]
\newtheorem{lemm}[theorem]{Lemma}
\newtheorem{coro}[theorem]{Corollary}
\newtheorem{prop}[theorem]{Proposition}
\newtheorem{namet}[theorem]{\myThmName}

\theoremstyle{remark}

\newtheorem{exas}[theorem]{Examples}
\newtheorem{rema}[theorem]{Remark}

\theoremstyle{definition}
\newtheorem{defs}[theorem]{Definitions}

\newcommand{\coloneqq}{\mathrel{\mathop:}=}
\newcommand{\id}{\mathrm{id}}
\newcommand{\smallset}[2]{\{{#1}\left\vert\vphantom{}\right.\,{#2}\}}
\newcommand{\Aut}[1]{\mathrm{Aut}(#1)}
\newcommand{\Cs}[2][]{\mathrm{C}_{#1}(#2)}
\newcommand{\Cg}[1]{\mathrm{C}_{#1}}
\newcommand{\trs}[2][]{\Gamma\ifx\empty#1\else^{[#1]}\fi_{[#2]}}
\newcommand{\trsp}[1]{\Gamma^{[#1]}}
\newcommand{\Fig}[1]{\operatorname{Fig}(#1)}
\newcommand{\DD}{\mathbb D}
\newcommand{\FF}{\mathbb F}
\newcommand{\HH}{\mathbb H}
\newcommand{\NN}{\mathbb N}
\newcommand{\UU}{\mathbb U}
\newcommand{\cB}{\mathcal{B}}
\newcommand{\cD}{\mathcal{D}}
\newcommand{\cH}{\mathcal{H}}
\newcommand{\cL}{\mathcal{L}}
\newcommand{\SL}[2][]{\mathrm{SL}_{#1}{(#2)}}
\newcommand{\PSL}[2][]{\mathrm{PSL}_{#1}{(#2)}}
\newcommand{\PgL}[2][]{\mathrm{P}\Gamma\mathrm{L}_{#1}{(#2)}}
\newcommand{\SU}[2][]{\mathrm{SU}_{#1}(#2)}
\newcommand{\PSU}[2][]{\mathrm{PSU}_{#1}(#2)}
\newcommand{\PG}[2]{\mathrm{PG}(#1,#2)}
\newcommand{\AG}[2]{\mathrm{AG}(#1,#2)}
\newcommand{\Sz}[1]{\mathrm{Sz}{(#1)}}

\begin{document}

\title{Unitals with many involutory translations}
\author{Theo Grundh\"ofer, Markus J. Stroppel, Hendrik {Van Maldeghem}}

\makeatletter
\newcommand{\subjclass}[2][1991]{%
  \let\@oldtitle\@title%
  \gdef\@title{\@oldtitle\footnotetext{#1 \emph{Mathematics subject classification.} #2}}%
}
\newcommand{\keywords}[1]{%
  \let\@@oldtitle\@title%
  \gdef\@title{\@@oldtitle\footnotetext{\emph{Key words and phrases.} #1.}}%
}
\makeatother




\date{}

\maketitle

\abstract{%
  \noindent%
  If every point of a unital is fixed by a non-trivial translation and
  at least one translation has order two then the unital is classical
  (i.e., hermitian).
}

\section{Introduction}\label{sec1}

A unital $\UU=(U,\cB)$ of order~$q$ is a $2$-$(q^3+1, q+1, 1)$-design;
i.e., an incidence structure with $\vert U\vert =q^3+1$ such that every block
$B\in\cB$ has exactly $q+1$ points, and any two points in~$U$ are
joined by a unique block in~$\cB$. It follows that every point is
incident with exactly~$q^2$ blocks. Without loss of generality, we
will assume throughout that the blocks are subsets of~$U$. %
The classical examples of unitals are the hermitian ones: For any
prime power~$q$, the points of the hermitian unital
$\cH(\FF_{q^2}\vert \FF_{q})$ are the absolute points with respect to a
suitable polarity of the projective plane~$\PG2{\FF_{q^2}}$ over the
field of order~$q^2$, see~\cite[II.8, pp.\,57--63]{MR0333959}
or~\cite[Ch.\,2]{MR2440325}.  The blocks are induced by secant lines
in~$\PG2{\FF_{q^2}}$.

A translation (with center~$c$) of a unital~$\UU$ is an automorphism
of~$\UU$ that fixes the point~$c$ and every block through~$c$.  We
write $\trs{c}$ for the group of all translations with center~$c$. %
Much is known for the case where there are at least two points $u,v$
such that $\trs{u}$ and $\trs{v}$ both have order~$q$,
see~\cite{GrundhoeferStroppelVanMaldeghem-moufangsets}. If there even
exist three non-collinear points with that property then the unital is
the classical hermitian one, see~\cite{MR3090721}.  We show
(Theorem~\ref{allCenters}) that the same conclusion holds if $\trs{c}$
has even order for every point~$c$. %
In fact, it suffices that every point is the center of some
non-trivial translation, and that at least one of these translations
is an involution.

Polar unitals in Figueroa planes of even order provide examples where
the centers of involutory translations form a proper subset not
contained in a block (namely, a subunital that is in fact hermitian),
see Theorem~\ref{translationsFigueroa}. %
Examples where the set of centers is contained in a block exist in
abundance, see Section~\ref{sec:fewTranslations}.

\newpage
\section{Transitivity}

The following elementary observation can be traced back to Gleason;
cp.~\cite[pp.\,190\,f]{MR0233275}, see~\cite[4.1]{MR3090721} for a
(very short) proof.

\begin{lemm}\label{Gleason}
  Let\/~$p$ be a prime, and let\/~$H$ be a group acting on a finite
  set\/~$X$.  Assume that for every $x\in X$ there exists in~$H$ an
  element of order~$p$ which fixes~$x$ but no other element of\/~$X$.
  Then\/~$H$ is transitive on~$X$. %
  \qed
\end{lemm}

Now let $\UU = (U,\cB)$ be a unital.  A translation of order~$n$ is
called an $n$-translation, and $\trs[n]c$ denotes the set of all
$n$-translations of~$\UU$ with center~$c$. %
For each~$n$, let\/~$\Omega_n$ be the (possibly empty) set of all
centers of $n$-translations, and
$\trsp{n} \coloneqq
\bigl\langle\bigcup_{c\in\Omega_n}\trs[n]c\bigr\rangle$.

Recall from~\cite[Theorem\,1.3]{MR3090721} that a non-trivial
translation of a unital fixes no point apart from its center. %
If~$\Omega_n$ is not empty, the size of any $\Aut\UU$-orbit
in~$\Omega_n$ is therefore congruent~$1$ modulo~$n$ (because a cyclic
group of size~$n$ acts semi-regularly on the complement of a point
in~$\Omega_n$), and every $\Aut\UU$-orbit outside~$\Omega_n$ has size
divisible by~$n$ (because the same cyclic group acts semi-regularly on
that set).

\begin{lemm}\label{trsOmegaP}
  Let\/ $p$ be a prime. %
  If\/ $H \le \Aut\UU$ contains a $p$-translation then~$H$ is
  transitive on the set of all centers of $p$-translations in~$H$. %
  In particular, we have:
  \begin{enumerate}
  \item\label{trsOmega}%
    The group $\trsp{p}$ %
    generated by all $p$-translations is transitive on~$\Omega_p$.
  \item\label{blockstabtrs}%
    For each block\/~$B$ of\/~$\UU$ the group %
    generated by all $p$-translations with center on~$B$ is transitive
    on the set\/~$\Omega_p\cap B$ of centers of $p$-translations
    in~$B$. %
    Consequently, the stabilizer 
    of~$B$ in $\trsp{p} \le \Aut\UU$ is transitive
    on\/~$\Omega_p\cap B$. %
  \item\label{composite}%
    If~$n>1$ is an integer such that\/~$\Omega_n$ is not empty, then
    the group $\trsp{n}$ is transitive on~$\Omega_n$ and\/
    $\Omega_k=\Omega_n$ holds for each divisor\/ $k>1$ of\/~$n$.
  \end{enumerate}
\end{lemm}
\begin{proof}
  The group~$H$ acts on the set of all centers of $p$-translations
  in~$H$, and Lemma~\ref{Gleason} applies. For assertions~\ref{trsOmega}
  and~\ref{blockstabtrs} we specialize $H=\trsp{p}$ and
  $H=\bigl\langle\bigcup_{x\in B}\trs[p]x\bigr\rangle$, respectively.
  Assertion~\ref{composite} follows from the fact that that~$\Omega_n$
  is contained in the $\Aut\UU$-orbit~$\Omega_k$.
\end{proof}

\begin{defs}
  Assume that $\UU$ admits a translation of order $k>1$. %
  Let $p_k$ be the smallest prime divisor of~$k$. Then
  $\Omega_k=\Omega_{p_k}$ holds by Lemma~\ref{trsOmegaP}.\ref{composite}. %
  Put $K \coloneqq \smallset{p_k}{k\in\NN,\vert \trsp{k}\vert >1}$.
  Clearly, we obtain a disjoint union %
  $U = \mho \cup \bigcup_{p\in K}\Omega_p$, where
  $\mho \coloneqq \smallset{x\in U}{\trs{x}=\{\id\}}$.
  
  We write $\cB_k$ for the set of all blocks of~$\UU$ that contain at
  least two points of~$\Omega_k$, and consider the incidence structure
  $\UU_k \coloneqq (\Omega_k,\cB_k,\in)$. If $\Omega_k$ is not
  contained in a block then~$\UU_k$ is a non-trivial linear space.  

  From Lemma~\ref{trsOmegaP}.\ref{trsOmega}
  and Lemma~\ref{trsOmegaP}.\ref{blockstabtrs} we know that~$\trsp{k}$
  induces a transitive group of automorphisms of~$\UU_k$, and that for
  each block $B\in\cB_k$ the stabilizer of~$B$ in~$\trsp{k}$ acts
  transitively on $\Omega_k\cap B$. %
\end{defs}

Recall that a substructure of an incidence structure (with points and
blocks) is called \emph{ideally embedded} if the pencils of blocks in
the substructure are the same as the pencils in the larger structure.

\begin{lemm}\label{idealEmb}
  Assume that\/ $\UU$ is a unital such that every point is the center
  of some non-trivial translation.  For each $p\in K$, the linear
  space~$\UU_p$ is then ideally embedded in~$\UU$. %
  In particular, the set $\Omega_p$ is not contained in a block.  The
  action of~$\Aut\UU$ on~$\Omega_p$ is faithful.

  If there exists $p\in K$ such that\/ $\vert \Omega_p\cap B\vert $ is constant
  for~$B\in\cB_p$ then $\Omega_p=U$, and\/ $\trsp{p}$ is transitive
  on~$U$. %
  In particular, this happens if\/~$\Aut\UU$ is two-transitive
  on~$\Omega_p$.
\end{lemm}
\begin{proof}
  Consider $p\in K$, and a block $B\in\cB$ through a point
  $x\in\Omega_p$. It suffices to consider the case where there exists
  a point $y\in B\smallsetminus\Omega_p$. By our assumption, there exists a
  non-trivial translation~$\tau$ with center~$y$. Now $x^\tau$ is a
  point different from~$x$, and lies in $\Omega_p\cap B$. %
  So $\UU_p$ is ideally embedded in~$\UU$. %
  Clearly, this implies that~$\Omega_p$ is not contained in any block.
    
  Assume that $\alpha\in\Aut\UU$ fixes every point
  in~$\Omega_p$. Then~$\alpha$ fixes each block that meets~$\Omega_p$
  because~$\UU_p$ is ideally embedded. Each point $x$
  outside~$\Omega_p$ lies on more than one block meeting~$\Omega_p$,
  so~$\alpha$ fixes every point of~$\UU$, and is trivial.

  Now assume that $b \coloneqq \vert \Omega_p\cap B\vert $ is constant
  for~$B\in\cB_p$; then $\vert \Omega_p\vert  = 1+q^2b$.  If there exists
  $z\in U\smallsetminus\Omega_p$ then, by our assumption, there is a
  translation of prime order~$r$ and center~$z$.  %
  Now $r$ divides both~$q$ and~$\vert \Omega_p\vert  = 1+q^2b$, and we obtain a
  contradiction.
\end{proof}

The Figueroa unitals (see Theorem~\ref{theo:Figueroa}) of even order are
examples that show that the condition $\mho=\emptyset$ is necessary
in Lemma~\ref{idealEmb}: in those unitals $\UU_2$ is not ideally embedded,
and there exists an automorphism of order~$3$ that acts trivially
on~$\Omega_2$. See also Proposition~\ref{subunital}.

\enlargethispage{5mm}%
\section{Results from group theory}

\begin{lemm}\label{sharplytrs}
  Let $G$ be a transitive permutation group on some finite
  set~$\Omega$ with more than one element. Suppose that the
  stabilizer~$G_x$ of some $x\in\Omega$ contains an involution~$\tau$
  that is semi-regular on $\Omega\smallsetminus\{x\}$, %
  and that\/~$G$ has a transitive normal subgroup~$M$ of odd order. %
  Then the following are equivalent:
  \begin{enumerate}
  \item $M$ is abelian.
  \item $M$ acts regularly (i.e., sharply transitively) on~$\Omega$.
  \item The involution~$\tau$ acts semi-regularly on~$M$ by
    conjugation.
  \end{enumerate}
  If one (and then any) of these conditions is satisfied then
  $\langle\tau\rangle M$ is a generalized dihedral group; i.e.,
  conjugation by~$\tau$ inverts each element of~$M$, and
  $\langle\tau\rangle M = \tau^M\cup M$.
\end{lemm}
\begin{proof}
  If~$M$ is abelian then the stabilizer~$M_x$ also fixes each element
  of the orbit~$\Omega$ of~$x$ under~$M$.  Hence~$M_x$ is trivial, and
  the action of~$M$ is regular. %
  If~$M$ acts regularly, we identify $m\in M$ with the
  image~$x^m$. Semi-regularity of~$\tau$ on~$\Omega$ then translates
  into semi-regularity of the automorphism induced by conjugation
  with~$\tau$ on~$M$. %
  The map $m\mapsto m^{-1}m^\tau$ is then injective, and bijective
  since~$M$ is finite. Write $m\in M$ as $m=k^{-1}k^\tau$ with
  $k\in M$, and then calculate
  $m^{-1} = (k^\tau)^{-1} k = (k^{-1})^\tau k = (k^{-1}k^\tau)^\tau =
  m^\tau$. %
  Thus the automorphism induced by~$\tau$ is the anti-automorphism
  $m\mapsto m^{-1}$, and~$M$ is abelian.

  If~$\tau$ induces inversion on~$M$ then
  $m^{-1}\tau m = \tau(\tau m^{-1}\tau)m = \tau m^2$ holds for each
  $m\in M$.  As $M$ has odd order, this means that the coset~$\tau M$
  equals the conjugacy class~$\tau^M$, and  consists of involutions.
\end{proof}

The following result has been  proved in
\cite[Theorem\,2]{MR322020}; cp. \cite[Theorem\,2]{MR318308}.

\begin{theorem}
  \label{hering}
  Assume that the group~$G$ acts transitively (but not necessarily
  faithfully) on some finite set~$\Omega$ with more than one
  element. Suppose that the stabilizer~$G_x$ of some $x\in\Omega$ has
  a normal subgroup~$Q$ of even order that is semi-regular on
  $\Omega\smallsetminus\{x\}$. %
  Then the normal closure $S \coloneqq \langle Q^G\rangle$ either has
  a transitive normal subgroup of odd order, or acts two-transitively
  on~$\Omega$ as one of the groups $\SL[2]{\FF_{2^e}}$,
  $\Sz{2^{2e-1}}$, or $\PSU[3]{\FF_{2^{2e}}\vert \FF_{2^e}}$ for
  suitable~$e\ge2$; the action is the usual two-transitive one. %

  The group $S$ itself is then $S=QN$ where $N$ denotes the largest
  normal subgroup of odd order in~$S$, or~$S$ is isomorphic to
  $\SL[2]{\FF_{2^e}}$, $\Sz{2^{2e-1}}$, $\SU[3]{\FF_{2^{2e}}\vert \FF_{2^e}}$ or
  $\PSU[3]{\FF_{2^{2e}}\vert \FF_{2^e}}$, according to the group induced
  on~$\Omega$. %
  \qed
\end{theorem}

We add some information contained in~\cite{MR322020} ---in particular,
see~\cite[Lemma\,3]{MR322020}--- but not in the statement of the theorem
referred to above. %

\begin{prop}\label{normalSubgroupAbelian}%
  The kernel of the action of~$G$ as considered in
  Theorem~\ref{hering} is $K=\Cs[G]S$. %

  If\/~$S$ has a transitive normal subgroup~$N$ of odd order then~$Q$
  acts fixed-point-freely on~$N/K$, so the group~$N/K$ induced by~$N$
  on~$\Omega$ is a sharply transitive abelian group in that case. %
  The group~$Q$ is then (isomorphic to) a Frobenius complement.

  If~$S$ induces a non-solvable group on~$\Omega$ then either
  $K\cap S$ is trivial, or $\vert K\cap S\vert =3$; the latter case can only
  occur if $S\cong\SU[3]{\FF_{2^{2e}}\vert \FF_{2^e}}$ with odd $e>1$.
  \qed
\end{prop}

The Sylow subgroups of a Frobenius complement are either cyclic or
generalized quaternion groups; see~\cite[10.3.1]{MR569209}. So
Proposition~\ref{normalSubgroupAbelian} yields:
  
\begin{coro}\label{onlyOneInvolution}
  If\/ $S= \langle Q^G\rangle$ has a transitive normal subgroup~$N$ of
  odd order then the group~$Q$ contains exactly one involution. %
  That involution induces inversion on~$N$. 
\end{coro}

\goodbreak%
\section{Unitals with involutory translations}

By Lemma~\ref{trsOmegaP}, Hering's result (see Theorem~\ref{hering})
applies to the permutation group induced by $\trsp2$ on~$\Omega_2$ if
this set has more than one element.

The following proposition is a small part of the 
results in~\cite{MR773556}. %
For the reader's convenience, we include a proof of the facts that we
need here. %

\begin{prop}\label{prop:linearSpaces}
  Let $(X,\cL)$ be a linear space with~$v$ points such that each line
  has~$k$ points. %
  Assume that\/ $k>2$. %
  Consider a group~$G$ isomorphic to one of the groups
  $\PSL[2]{\FF_q}$, $\Sz{q}$ or $\PSU[3]{\FF_{q^2}\vert \FF_q}$, for
  some prime power~$q$.  If\/~$G$ acts on~$X$ in its usual
  two-transitive action and by automorphisms of $(X,\cL)$, then either
  $v=k$ (and there is just one line in~$\cL$), or
  $G\cong\PSU[3]{\FF_{q^2}\vert \FF_q}$ and\/
  $(X,\cL) \cong\cH(\FF_{q^2}\vert \FF_q)$ is isomorphic to the
  hermitian unital of order~$q$.
\end{prop}

\begin{proof}
  We assume that there is more than one line, and discuss the three
  different cases separately.

  \textbf{(PSL) } %
  The usual two-transitive action of $G\cong\PSL[2]{\FF_{q}}$ is the
  natural action on the projective line $\FF_q\cup\{\infty\}$ via
  fractional linear transformations. Then $v=q+1$, and the stabilizer
  of two points $x,y\in X$ has two orbits of length~$1$ and at most
  two other orbits, each of length $\frac{q-1}{\gcd(2,q-1)}$. The line
  joining~$x$ and~$y$ contains at least one of those orbits, and
  $k-1 \ge \frac{q-1}{\gcd(2,q-1)}+1 > \frac{v-1}2$.  So there is no
  space left for a second line through~$x$, and $v=k$ follows.

  The Suzuki groups and the unitary groups need a closer look; we will
  use the following facts about linear spaces with $v$ points and
  $k>2$ points per line:
      
  The number of lines per point in $(X,\cL)$ is $r =
  \frac{v-1}{k-1}$. %
  We assume that the linear space has more than one line, so $r>1$. %
  Fisher's inequality (see~\cite[1.3.8, p.\,20]{MR0233275}) then says
  $r\ge k$, and equality holds precisely if the linear space is a
  projective plane. In the latter case, we have $v=k^2-k+1$. If $r>k$
  then $v-1 \ge (k+1)(k-1) = k^2-1$.

  \textbf{(Sz) } %
  The usual two-transitive action of $G \cong \Sz{q}$ is its action on
  the Suzuki-Tits ovoid; see~\cite[Section\,21]{MR572791}.  Then
  $v=q^2+1$ and $q$ is a power of~$2$.  Hence we cannot have
  $v-1=k^2-k$, and $v-1\ge k^2-1$ follows.

  For any two points $x,y\in X$, the stabilizer $G_{x,y}$ acts on~$X$
  with two orbits of length~$1$, every other orbit has length $q-1$;
  see~\cite[Lemma\,21.5]{MR572791}. %
  The block joining $x$ and~$y$ is thus the union of $\{x,y\}$ with a
  collection of such orbits of length $q-1$, and
  $k \ge (q-1)+2 = q+1$. %
  Now the inequality %
  $q^2 = v-1 \ge k^2-1 \ge (q+1)^2-1$ yields a contradiction.

  \textbf{(PSU) } %
  The usual two-transitive action of
  $G \cong \PSU[3]{\FF_{q^2}\vert \FF_{q}}$ is its action on the points of
  the hermitian unital of order~$q$; see~\cite{MR175980},
  \cite{MR0295934}, \cite[Theorem\,5.2]{MR0333959}.  Then $v=q^3+1$, the
  fact that $q$ is a prime power yields $v-1 \ne k^2-k$, and we have
  $q^3 = v-1 \ge k^2-1$.

  If $q=2$ then $v=9$ and $k>2$ together with $k^2-1 \le v-1$ yields
  $k=3$. Then $(X,\cL)$ is determined uniquely, we have
  $(X,\cL) \cong \cH(\FF_4\vert \FF_2) \cong \AG2{\FF_3}$. We assume $q>2$
  from now on.

  The stabilizer $G_{x,y}$ of two points $x,y\in X$ has two orbits of
  length~$1$, one orbit of length $q-1$, and every other orbit has
  length $\frac{q^2-1}{\gcd(3,q+1)}$;
  see~\cite[p.\,499]{MR0295934}. %
  So there are integers $s\in\{0,1\}$ and~$t\ge0$ such that
  $k-1 = s(q-1) + t \frac{q^2-1}{\gcd(3,q+1)} +1$.

  Aiming at a contradiction, we assume $t>0$. %
  Then $k \ge {\frac{q^2-1}{\gcd(3,q+1)} +2} \ge \frac{q^2+5}{3}$.
  From %
  $q^3 = v-1 \ge k^2-1 \ge \bigl(\frac{q^2+5}3\bigr)^2 -1$ we thus
  obtain $0 \ge q^4 -9q^3 + 10q^2 + 16 > q^2 \bigl(q(q-9)
  +10\bigr)$. %
  So $q<8$, and $q\in\{3,4,5,7\}$.  For $q\in \{3,4,7\}$ we have
  $\gcd(3,q+1)=1$, and our assumption $t>0$ yields
  $q^3 \ge \left((q^2-1)+1\right)^2 = q^4$. This is impossible.  For
  $q=5$ and $t>1$, we obtain the contradiction
  $5^3 \ge \bigl(2\frac{(5^2-1)}3+1\bigr)^2 = 17^2$.  So $t=1$, and
  $k-1 = 4s+8+1$ divides $v-1=5^3$. %
  Both cases for $s\in\{0,1\}$ lead to a contradiction.

  Therefore, we have $t=0$ and~$s=1$; and the block through $x$
  and~$y$ is the union of $\{x,y\}$ with the unique orbit of length
  $q-1$ under~$G_{x,y}$. This means that $(X,\cL)$ is isomorphic to
  the hermitian unital $\cH(\FF_{q^2}\vert \FF_{q})$.
\end{proof}

\begin{prop}\label{subunital}
  Let\/~$\UU$ be a unital of order~$q$, and assume that\/~$\Omega_2$
  is not contained in a block.  If the group induced by\/~$\trsp2$
  on~$\Omega_2$ does not have a normal subgroup that acts regularly
  on~$\Omega_2$ then that induced group is isomorphic
  to~$\PSU[3]{\FF_{2^{2e}}\vert \FF_{2^e}}$, and\/
  $\UU_2 = (\Omega_2,\cB_2)$ is isomorphic to the hermitian unital of
  order~$2^e$, with $e>1$. %
\end{prop}
\begin{proof}
  From Theorem~\ref{hering} we know that the group induced
  by\/~$\trsp2$ on~$\Omega_2$ is isomorphic to one of the groups
  $\SL[2]{\FF_{2^e}}$, $\Sz{2^{2e-1}}$, or
  $\PSU[3]{\FF_{2^{2e}}\vert \FF_{2^e}}$ (for suitable~$e\ge2$), with
  the usual two-transitive action. %
  From Proposition~\ref{prop:linearSpaces} we then know that
  $\trsp2 \cong \PSU[3]{\FF_{2^{2e}}\vert \FF_{2^e}}$, and~$\UU_2$ is
  the hermitian unital of order~$2^e$.
\end{proof}

The situation in Proposition~\ref{subunital} actually occurs, for
example, in polar unitals of Figueroa planes of even order, see
Theorem~\ref{theo:Figueroa} and Theorem~\ref{translationsFigueroa}
below.  Our Proposition~\ref{subunital} is a version for unitals (with
translations) of a result \cite[Theorem\,5.2]{MR0467511} on
projective planes (and elations).
  
\begin{theorem}\label{allCenters}%
  Let\/ $\UU=(U,\cB)$ be a unital of order~$q$, assume that every
  point of\/~$\UU$ is the center of some non-trivial translation, and
  that there exists a translation of order~$2$. %
  Then~$q$ is a power of\/~$2$, and the unital\/~$\UU$ is the hermitian
  unital of order~$q$. 
\end{theorem}
\begin{proof}
  If $q=2$ then the unital is the hermitian unital of order~$2$, see 
  \cite[10.16]{MR1189139}. %
  We assume $q>2$ in the rest of the proof.

  \textbf{Case A: } %
  Assume first that $\Omega_2=U$, i.e., every point of~$\UU$ is the
  center of some involutory translation. %
  If the group $\trsp2$ generated by the involutory translations does
  not have a regular normal subgroup then Proposition~\ref{subunital}
  says %
  that $\trsp2 \cong \PSU[3]{\FF_{2^{2e}}\vert \FF_{2^e}}$, and $\UU$
  is the hermitian unital of order~$2^e$.

  (Alternatively, this can be derived more directly from
  Theorem~\ref{hering}: the groups $\SL[2]{\FF_{2^e}}$ are excluded
  since they are triply transitive; the orbits of the
  two-point-stabilizers in $\Sz{2^{2e-1}}$ are too large to yield
  blocks of the unital. For the groups
  $\PSU[3]{\FF_{2^{2e}}\vert \FF_{2^e}}$ the block through two points
  consists of those two together with the unique shortest orbit of
  their stabilizer.  Hence $\UU$ coincides with the hermitian unital
  $\cH(\FF_{2^{2e}}\vert \FF_{2^e})$ of order~$2^e$.)
  
  If the group~$\trsp2$ has a regular normal subgroup~$N$, then~$N$ is
  abelian (see Proposition~\ref{normalSubgroupAbelian}).  Consider two blocks
  $B,C$ through some point. Then the stabilizers $N_B$ and~$N_C$ are
  subgroups of order $q+1$ in~$N$, and $N_BN_C$ is a subgroup of order
  $(q+1)^2$. So $(q+1)^2$ divides $\vert N\vert =q^3+1$, and $q+1$ divides
  $(q^3+1)/(q+1) = q^2-q+1 = (q+1)(q-1)-(q-2)$. This implies $q-2=0$,
  contradicting our assumption $q>2$ (in fact, $\trsp2$ on the unital
  of order~$2$ does contain a regular normal subgroup,
  see Remark~\ref{specialOrderTwo}).

  \textbf{Case B: } %
  Now we assume $\Omega_2\ne U$, and aim at a contradiction. %
  As every point is the center of a non-trivial translation, there is
  a prime~$p$ such that $\Omega_p$ is not empty, and disjoint
  to~$\Omega_2$. %
  We know from Lemma~\ref{idealEmb} that $\Omega_2$ is not contained in a
  block, and that the group~$\trsp2$ generated by all involutory
  translations acts faithfully on~$\Omega_2$. %
    
  The action of~$\trsp2$ on~$\Omega_2$ is as in Theorem~\ref{hering}
  and Proposition~\ref{normalSubgroupAbelian}. %
  If $\trsp2$ is two-transitive on~$\Omega_2$ then
  Lemma~\ref{idealEmb} yields $\Omega_2=U$, contradicting our present
  assumptions. %
  If $\trsp2$ is not two-transitive on~$\Omega_2$ then~$\trsp2$ has an
  abelian normal subgroup~$A$ acting sharply transitively
  on~$\Omega_2$, see Theorem~\ref{hering}. %
  From Corollary~\ref{onlyOneInvolution} we know that each one of the
  involutory translations acts by inversion on~$A$, and that $\trs{x}$
  contains exactly one involution if $x\in\Omega_2$. %
  We denote that involution by~$j_x$. %
  The order of~$A$ equals that of~$\Omega_2$, so it is odd, and
  divisible by~$p$.

  Consider $a\in A$ and an arbitrary point $x\in\Omega_2$. Then both
  $j_x$ and $a^{-1}j_xa = j_{x^a}$ fix the block joining~$x$
  and~$x^a$. So $a^2 = j_xa^{-1}j_xa$ fixes that block, and so
  does~$a$ because squaring is an automorphism of the abelian
  group~$\langle a\rangle$ of odd order.

  As $\vert A\vert =\vert \Omega_2\vert $ is divisible by~$p$, we find an element $a\in A$
  of order~$p$.  That~$a$ fixes a point $c\in\Omega_p$ because
  $\vert \Omega_p\vert \equiv1\pmod{p}$. %
  For any $x\in\Omega_2$, we have $a^2=j_xa^{-1}j_xa$, and $c=c^{a}$
  yields that $x^a$ lies on the block joining~$x$ and~$c$ because the
  translations~$j_x$ and~$a^{-1}j_xa$ fix each block through their
  respective center. %
  So~$a$ fixes every block joining a point of~$\Omega_2$ with~$c$.
  As~$a$ fixes no point in~$\Omega_2$, the point~$c$ is the only point
  fixed by~$a$.

  Let~$j$ be any involutory translation. 
  As $c$ is fixed by $a^2=ja^{-1}ja$, we obtain
  $c^{j} = c^{a^{-1}ja} = c^{ja}$, and the point $c^{j}$ is
  fixed by~$a$. Since~$c$ is the only fixed point of~$a$, we reach the
  contradiction that the involutory translation~$j$ fixes
  $c\notin\Omega_2$. 
\end{proof}

\begin{rema}\label{specialOrderTwo}
  Let $p$ be a prime, and let $e$ be a positive integer. Then each
  non-trivial translation of the hermitian
  unital~$\cH(\FF_{p^{2e}}\vert \FF_{p^e})$ of order~$p^e$ has order~$p$.
  If $p^e>2$ then the group~$\trsp{p}$ generated by all translations
  is simple, coincides with $\PSU[3]{\FF_{p^{2e}}\vert \FF_{p^e}}$, and
  acts two-transitively on the point set
  of~$\cH(\FF_{p^{2e}}\vert \FF_{p^e})$; see~\cite[10.15, 10.12]{MR1189139}.

  On the hermitian unital of order~$2$, the group $\trsp2$ behaves in
  an exceptional way: %
  then the group~$\trsp2$ is solvable, and not two-transitive, but it
  is still the commutator subgroup of $\PSU[3]{\FF_{4}\vert \FF_{2}}$;
  see~\cite[10.15 and the discussion on pp.\,123f]{MR1189139}. %
    
  The smallest unital (of order~$2$) is isomorphic to the affine plane
  of order~$3$. What we call a translation of the unital is an affine
  homology in that plane; this explains the structure
  of~$\trsp2 \cong \Cg2\ltimes\Cg3^2$, a generalized dihedral group. %
  That group is transitive on the point set, but not transitive on the
  block set (it preserves each parallel class in the affine plane).
\end{rema}

\section{Figueroa unitals}

\begin{exas}\label{exas:Figueroa}
  Let $r$ be a prime power. \cite{MR682668} has constructed a
  projective plane $\Fig{r^3}$ of order~$r^3$ with a pappian
  subplane~$\DD$ of order~$r$; see~\cite{MR692241} for the case of
  general~$r$, and~\cite{MR850165} for a synthetic construction. %
  The plane~$\Fig{r^3}$ is not desarguesian unless $r=2$.  The full
  group of automorphisms of the subplane~$\DD$ extends to a group
  $G \cong \PgL[3]{\FF_{r}}$ of automorphisms of~$\Fig{r^3}$.  If
  $r>2$ then the full automorphism group $\Aut{\Fig{r^3}}$ is the
  direct product of a group of order~$3$ with said group~$G$;
  see~\cite{MR692241}, cp.~\cite{MR818716}. The cyclic factor is
  generated by a planar automorphism~$\alpha$ (of order~$3$) that is
  used in the construction of~$\Fig{r^3}$. %
  Dempwolff~\cite[Theorem\,B]{MR787308} has noted that every elation
  of the subplane~$\DD$ (as an element of the group~$G$) is induced by
  an elation of~$\Fig{r^3}$.

  \goodbreak%
  If $r$ is a perfect square, say $r=q^2$, then there is a
  polarity~$\pi$ of the Figueroa plane of order~$q^6$, and the
  absolute points of~$\pi$ carry a unital~$\UU_{\Fig{q^6}}$ of order
  $q^3$, see~\cite{MR1607945}. %
  The unital $\UU_{\Fig{q^6}}$ is not hermitian;
  see~\cite{MR2995129}.  In fact, there are O'Nan configurations,
  see~\cite{MR3208081}.  The intersection~$\HH$ of~$\UU_{\Fig{q^6}}$
  with the subplane~$\DD$ is isomorphic to the hermitian unital of
  order~$q$.

  The centralizer of the polarity~$\pi$ in~$\Aut{\Fig{q^6}}$ is the
  direct product of $\langle\alpha\rangle$ with the centralizer
  of~$\pi$ in $G \cong \PgL[3]{\FF_{q^2}}$. %
  In particular, every translation of~$\HH$ is induced by an elation
  of the subplane~$\DD$, and thus by an elation~$\psi$ of~$\Fig{q^6}$
  in the centralizer of~$\pi$. The restriction of~$\psi$
  to~$\UU_{\Fig{q^6}}$ is then a translation of~$\UU_{\Fig{q^6}}$.
\end{exas}

We obtain: %
  
\enlargethispage{5mm}%
\begin{theorem}\label{theo:Figueroa}
  In the unital\/~$\UU_{\Fig{q^6}}$ of order~$q^3$, there is a
  hermitian subunital\/~$\HH$ of order~$q$ such that every point
  of\/~$\HH$ is the center of a group of order~$q$ consisting of
  translations. %
  In particular, the point set of\/~$\HH$ is contained in~$\Omega_p$.

  Let\/ $G$ be the group generated by all elations of the Figueroa
  plane $\Fig{q^6}$ that leave $\UU_{\Fig{q^6}}$ invariant. %
  Then~$G$ is isomorphic to the commutator group %
  of\/ $\PSU[3]{\FF_{q^2}\vert \FF_{q}}$, and acts two-transitively on the
  point set of\/~$\HH$.
  \qed
\end{theorem}

\begin{theorem}\label{translationsFigueroa}
  Let\/ $q$ be a power of~$2$. Then every involutory translation
  of\/~$\UU_{\Fig{q^6}}$ has its center in~$\HH$, and is induced by an
  elation of the Figueroa plane\/~$\Fig{q^6}$ that
  leaves~$\UU_{\Fig{q^6}}$ invariant. %
  In particular, we have $\HH = \UU_2 = (\Omega_2,\cB_2)$, the
  subunital~$\HH$ is invariant under $\Aut{\UU_{\Fig{q^6}}}$, and
  every non-trivial translation is an involution. %
\end{theorem}
\begin{proof}
  We write $q=2^{a}$ with $a\in\NN$; then $\UU_{\Fig{q^6}}$ has
  order~$2^{3a}$ and~$\HH$ is a unital of order~$2^a$.  The order of
  any translation of $\UU_{\Fig{q^6}}$ divides~$q$, and is thus a
  power of~$2$. If the translation is not trivial then
  Lemma~\ref{trsOmegaP} yields that its center lies in~$\Omega_2$.
  
  From Theorem~\ref{theo:Figueroa} we know that the point set of~$\HH$
  is contained in~$\Omega_2$. %
  Assume that~$\Omega_2$ is not contained in the
  subunital\/~$\HH$. Then Proposition~\ref{subunital} yields
  that~$\trsp2$ induces a group isomorphic to
  $\PSU[3]{\FF_{2^{2e}}\vert \FF_{2^e}}$ on~$\Omega_2$, and $\UU_2$ is
  isomorphic to the hermitian unital of order~$2^e$.  Thus
  $2^a<2^e\le 2^{3a}$, and
  $\UU_2\cong\cH(\FF_{2^{6a}}\vert \FF_{2^{3a}})$ contains the
  subunital $\HH\cong\cH(\FF_{2^{2a}}\vert \FF_{2^{a}})$. According
  to~\cite{GrundhoeferStroppelVanMaldeghem-subunitals-arxiv}, the embedding
  of the unitals is given by an embedding of quadratic field
  extensions. This leaves only the possibility $e=3a$, but then
  $\UU_{\Fig{q^6}}=\UU_2$ is hermitian, a contradiction. %
  So $\UU_2=\HH$, every translation of $\UU_{\Fig{q^6}}$ has its
  center in~$\HH$, and~$\HH$ is invariant under all translations of
  $\UU_{\Fig{q^6}}$.

  If~$\tau$ is any non-trivial translation of $\UU_{\Fig{q^6}}$ then
  the order of $\tau$ divides~$q^3=2^{3a}$. So the center~$c$
  of~$\tau$ lies in~$\Omega_2$, and is thus a point of~$\HH$. %
  The elations of~$\Fig{q^6}$ with center~$c$ that
  leave~$\UU_{\Fig{q^6}}$ invariant form a group~$\Theta$ that acts
  faithfully on~$\UU_{\Fig{q^6}}$ and induces the full group of
  translations of~$\HH$ with center~$c$.

  Let $B$ be a block through~$c$, and consider any point
  $x\in(\Omega_2\cap B)\smallsetminus\{c\}$. %
  Since the group of all translations of the hermitian unital~$\HH$ is
  transitive on $(\Omega_2\cap B)\smallsetminus\{c\}$, there exists an
  involutory elation $\vartheta\in\Theta$ with
  $x^\tau=x^\vartheta$. Then $\tau\vartheta$ is a translation
  of~$\UU_{\Fig{q^6}}$ with center~$c$ fixing $x\ne c$, so
  $\tau\vartheta$ is trivial on~$\UU_{\Fig{q^6}}$, and
  $\tau=\vartheta$.
\end{proof}

\section{Examples: unitals with few translations}
\label{sec:fewTranslations}

There are unitals with no translations at all, e.g., the Ree unitals,
or the presently known unitals of order~$6$. %
See~\cite[1.8]{MR3090721} and~\cite{MR2795696} for the Ree unitals;
the unitals of order~$6$ are treated in~\cite[5.1,
p.\,2888]{MR3533336}; the information about the automorphisms given
there suffices to see that there are no translations. %
Most of the unitals of order~$3$ and many unitals of order~$4$ found
by computer do not admit automorphisms of order~$3$ or~$2$,
respectively, let alone translations; see~\cite{MR3138681},
\cite[Table\,III, p.\,301]{MR2838908}.
  
There are unitals of prime power order~$p^e$ where $\Omega_p$ consists
of a single point (and, obviously, $\Omega_r$ is empty for every prime
$r\ne p$). %
For instance, consider a Coulter-Matthews plane of order~$3^e$ with
even~$e$, defined by a suitable planar monomial;
see~\cite{MR1432296}. Such a plane admits a unitary polarity; the
absolute points carry a unital of order~$3^{e/2}$,
see~\cite[5.2]{MR2593326}. That unital has a point~$\infty$ with a
translation group~$\trs\infty$ of order~$3^e$,
see~\cite[6.2]{MR2593326}. These groups are elementary abelian
$3$-groups.

If the unital is not classical (this surely happens if a certain Baer
subplane is not desarguesian, see~\cite[6.8]{MR2593326}) then the
point~$\infty$ is fixed by every automorphism of the unital, and
$\infty$ is the only center of any translation (for abstract
automorphisms, this follows from a deep result~\cite{MR3090721}, it
does not suffice to note that the point~$\infty$ is fixed by every
automorphism of the plane).

In planes over finite Dickson semifields, and in planes over twisted
fields, one also finds non-classical unitals (polar and otherwise)
with exactly one center of translations, see~\cite{MR3133742}
and~\cite[5.2]{MR3533345}, respectively.

For each order~$q=p^d$ with~$p$ prime, M\"ohler gives a construction
of unitals \cite[4.1]{moehler2021affine} depending on the choice of a
suitable family~$\cD$ of subsets of~$\SL[2]{\FF_q}$ such that either
$\Omega_p$ is a block and $\trsp{p}\cong\SL[2]{\FF_q}$, or the unital
is hermitian; see~\cite[3.11]{moehler2020automorphisms-IIG}. Suitable
families~$\cD$ leading to non-hermitian unitals are known for
$q\in\{4,8\}$, see~\cite[2.1, 3.5]{MR3533345}
and~\cite[Sect.\,3]{moehler2021affine}.
In~\cite[Corollary\,3.8]{moehler2020parallelisms-JGEOM} it is proved:
For $q$ even, the set~$\Omega_2$ in Gr\"uning's unital
(\cite{MR895943}, see~\cite[5.5]{MR3533345} for the description
needed here) has size~$q+1$, and that unital admits exactly $q+1$
non-trivial translations, each of order~$2$. %
In~\cite[4.6, 4.7, 4.8]{moehler2020parallelisms-JGEOM} one finds
unitals of order~$4$ with no translations, unitals of order~$4$ with
$\vert \Omega_2\vert =1$ and $\trsp2\cong\Cg2$, and unitals of
order~$4$ with $\vert \Omega_2\vert =1$ and
$\trsp2\cong\Cg2\times\Cg2$.

\goodbreak%

\providecommand{\noopsort}[1]{}\def\cprime{$'$}
  \def\polhk#1{\setbox0=\hbox{#1}{\ooalign{\hidewidth
  \lower1.5ex\hbox{`}\hidewidth\crcr\unhbox0}}}

\end{document}